\magnification 1200
\baselineskip 14pt
\parskip=3pt plus1pt minus.5pt
%\hsize 17truecm

\def\qed{\hfill\vbox{\hrule\hbox{\vrule\kern3pt
\vbox{\kern6pt}\kern3pt\vrule}\hrule}}

\def\rationalmap{\hbox{-\hskip -.1mm -\hskip -.1mm -\hskip -.1mm -\hskip
-.1mm -\hskip -1.2mm}_{^>}}

\def\mapright#1{\smash{
   \mathop{\longrightarrow}\limits^{#1}}}

\def\P{{\bf P}}
\def\O{{\cal O}}
\def\L{{\Lambda}}
\def\l{\ell}
\def\C{{k}}

\def\Ad{[{\AA}d]}
\def\Yo{[Ar]}
\def\ABT{[ABT]}
\def\Arrondo-Sols{[A-S]}
\def\Fantechi{[F]}
\def\HR{[H-R]}
\def\Severi{[S]}
\def\Zakpaper{[Z1]}
\def\Zakbook{[Z2]}

\centerline{\bf PROJECTIONS OF GRASSMANNIANS OF LINES AND}
\centerline{\bf CHARACTERIZATION OF VERONESE VARIETIES}
\bigskip
\centerline{\it by Enrique Arrondo}
\centerline{Departamento de Algebra}
\centerline{Facultad de Ciencias Matem\'aticas}
\centerline{Universidad Complutense de Madrid}
\centerline{28040 Madrid, Spain}
\centerline{Fax number: +34-1-394.46.07}
\centerline{E-mail: {\tt enrique@sunal1.mat.ucm.es}}

\bigskip
\bigskip
ABSTRACT: We characterize the double Veronese embedding of ${\bf P}^n$ as the
only variety that, under certain general conditions, can be isomorphically
projected from the Grassmannian of lines in ${\bf P}^{2n+1}$ to the
Grassmannian of lines in ${\bf P}^{n+1}$.
\bigskip
\bigskip

{\bf 0. Introduction.}
\bigskip
In \Severi, Severi proved that the only nondegenerate (i.e. not contained
in a hyperplane) smooth (complex) surface in $\P^5$ that can be isomorphically
projected to $\P^4$ is the Veronese surface. More recently, Zak
extended Severi's result and proved that, for $n\ge2$ the only nondegenerate
$n$-dimensional smooth subvariety of $\P^{{n(n+3)\over2}}$ that
can be isomorphically projected to
$\P^{2n}$ is the $n$-uple Veronese embedding of $\P^n$ (see \Zakbook, or \Ad\
for a similar statement).

In \Arrondo-Sols, there is a classification of all smooth surfaces in
$G(1,3)$ (the Grassmann variety of lines in $\P^3$) that are
non-trivial projection of a surface in $G(1,4)$ (by {\it non-trivial} we
mean that the corresponding surface in $G(1,4)$ is {\it nondegenerate} in the
sense that there is no hyperplane in $\P^4$ containing all the lines
parametrized by the surface). In particular, this classification shows that
the only nondegenerate smooth surface in
$G(1,5)$ that can be isomorphically projected to $G(1,3)$ is a Veronese
surface, more precisely the embedding of $\P^2$ in $G(1,5)$ by the vector
bundle
$\O_{\P^2}(1)\oplus\O_{\P^2}(1)$.

More generally, consider the embedding of $\P^n$ in $G(1,2n+1)$ given by
the vector bundle $\O_{\P^n}(1)\oplus\O_{\P^n}(1)$ (in this context we will
refer to it as the {\it $n$-dimensional Veronese variety} in $G(1,2n+1)$ ).
We showed in \Yo\ (see also Example 1.1) that this Veronese variety can be
isomorphically projected into $G(1,n+1)$, and made the following:

\noindent {\bf Conjecture 0.1:} For any $n\ge1$, the only nondegenerate 
smooth complex $n$-dimensional subvariety $X$ of $G(1,2n+1)$ that can be
isomorphically projected into $G(1,n+1)$ is the Veronese variety.

In this paper we prove this conjecture (see Theorem 3.1) under the extra
assumption that
$X$ is what we will call {\it uncompressed}, i.e. that the union in 
$\P^{2n+1}$ of all lines parametrized by $X$ has the expected dimension $n+1$.
I have not been able to remove this condition, although I am sure that it is
not necessary. In fact, the conjecture is known to be true for $n\le 2$, and
we can also prove it in case $n=3$ (see Corollary 4.1). However the kind of
proof required for the compressed case seems to be completely different from
the techniques introduced  in this paper.

The steps in the proof will follow the same as in \Zakbook\ for the projective 
case. Surprisingly, the difficulty for each step in the Grassmannian case
seems to be complementary to the difficulty in the projective case. For
example, the most tricky part (probably the only one) in our case is to prove
that the projectability of a variety implies that the appropriate secant
variety has small dimension (Lemma 2.2). Our approach to this result consists
of an infinitesimal study, which makes our result to depend strongly on the
characteristic zero assumption. But on the other hand, the main point in Zak's
proof is the so-called Terracini's lemma, used to prove some tangency result.
However in our case this tangency condition (Lemma 2.5) follows immediately
from the geometry of the Grassmannian of lines. 

In some sense, Grassmannians of lines seem to provide a much more natural
context to study these projection properties. For instance, any $G(1,n+1)$ has
dimension $2n$, so it is natural to expect that few smooth $n$-dimensional
subvarieties of it are projected from bigger Grassmannians. Also our result
works even for the case of curves ($n=1$), in which the theorem of Zak does
not give a characterization of the corresponding Veronese variety (i.e. a
conic). 

The organization of the paper is as follows. In section 1 we give some
preliminaries and recall some facts from \Yo. In section 2 we prove some
lemmas we will need to prove our theorem. In section 3 we state and prove our
main theorem. Finally in section 4 we discuss the extra hypothesis we added 
to our theorem. We also discuss some general results of what could be a deeper
study of projection properties in any Grassmannian of lines. In fact I hope 
that this paper will be just the starting point for such a study, even for
general Grassmannians. 

I would like to acknowledge the support of the Spanish CICYT through the
grant PB 93-0440-C03-01. Also I want to thank Edoardo Ballico for his
hospitality when I visited Trento in March and April of 1996. There he 
encouraged me to give a seminar on subvarieties of Grassmannians, as well as
to write some notes on the subject (see \Yo). This provided me the convenient
atmosphere to start thinking of this problem. Finally I want to thank Fyodor
Zak for providing me several references and comments; he also pointed me
that a previous proof of Lemma 2.3 was not completely correct.

\bigskip
\bigskip

{\bf 1. Preliminaries.}
\bigskip

We will work over an algebraically closed field $k$ of characteristic zero. 
We will denote by $G(r,m)$ the Grassmann variety of $r$-linear spaces in
$\P^m$. A linear projection $G(1,m)\rationalmap G(1,m')$ will mean 
the natural (rational) map induced by the corresponding linear projection
$\P^m\rationalmap\P^{m'}$.

\noindent {\bf Notation:} We will denote the elements of any Grassmannian of 
lines by small letters, say $\l$, and use the corresponding capital letter, 
say $L$, for the line in projective space that they define.

\noindent {\bf Example 1.1:} For the sake of completeness, let us recall 
here the example provided by Proposition 3.4  in \Yo. Consider the natural
embedding of $\P^n$ in $G(1,2n+1)$ defined by
$\O_{\P^n}(1)\oplus\O_{\P^n}(1)$. In coordinates, it can be described by
associating to each $(t_0:\ldots:t_n)\in\P^n$ the line spanned by the rows of
the matrix
$$\pmatrix{t_0&\ldots&t_n&0&\ldots&0\cr
           0&\ldots&0&t_0&\ldots&t_n}$$
\noindent We consider now the linear projection $\P^{2n+1}\rationalmap
\P^{n+1}$ defined by
$$(x_0:\ldots:x_{2n+1})\mapsto (x_0:x_1+x_{n+1}:\ldots:x_n+x_{2n}:x_{2n+1})$$
\noindent This projection induces a projection from $G(1,2n+1)$ to $G(1,n+1)$
and the image of $\P^n$ corresponds to the lines spanned by the rows of the
matrix
$$\pmatrix{t_0&t_1&\ldots&t_n&0\cr
           0&t_0&\ldots&t_{n-1}&t_n}$$
This is still an embedding of $\P^n$ in $G(1,n+1)$, since the maximal minors 
of the above matrix (which give the image of $\P^n$ after the Pl\"ucker
embedding of $G(1,n+1)$ in
$\P^{{n+2\choose2}-1}$) define the double Veronese  embedding of $\P^n$ in
$\P^{{n+2\choose2}-1}$. By this reason we will call this subvariety of
$G(1,2n+1)$ (or of $G(1,n+1)$ ) the {\it $n$-dimensional Veronese variety}.

\bigskip

For the rest of the paper (except for the last section) our setting will be 
the one needed for proving our main theorem. Most of the results could be
formulated in a more general setting, but we will leave this kind of comments
for the last section. So we will consider $X$ to be a smooth irreducible
$n$-dimensional subvariety of $G(1,2n+1)$. For short we will sometimes call a 
{\it line of $X$} to a line in $\P^{2n+1}$ parametrized by a point of $X$. We
will say that $X$ is {\it nondegenerate} if the union of all the lines of $X$
is not contained in a hyperplane. We will also say that $X$ is {\it
compressed} if the union in $\P^{2n+1}$ of all of its lines has dimension at
most $n$. Otherwise, if this union has dimension is $n+1$, we will say that
$X$ is {\it uncompressed}.

\noindent{\bf Remark 1.2:} The notion of compressedness is related with the
following  fact. Consider the incidence variety $I_X:=\{(\l,p)\in
X\times\P^{2n+1}\ | \ p\in L\}$ and the projection $p_X:I_X\to \P^{2n+1}$.
Then $X$ is uncompressed if and only if this projection is generically finite
over its image. In other words, for a general hyperplane $H$ of $\P^{2n+1}$
the rational map
$p_H:X\rationalmap H$ that assigns to each line of $X$ (not contained in $H$)
its intersection with $H$ is generically finite over its image.

\bigskip
\bigskip

{\bf 2. Some previous results.}
\bigskip

{}From now on, $X$ will be a smooth irreducible $n$-dimensional
uncompressed nondegenerate subvariety of $G(1,2n+1)$ that can be isomorphically
projected to $G(1,n+1)$. This means in particular that there is an
$(n-1)$-dimensional center of projection $\L$ verifying the following two
conditions:

\item{(*)} Any two skew lines of $X$ (probably infinitely close) span a
three-dimensional linear space that meets $\L$ at
most in one point.
\item{(**)} If two lines of $X$ (probably infinitely close) span only a
two-dimensional linear space, then this span does not meet $\L$.

Condition (*) is easier to handle in the sense that can be described in
terms of an irreducible variety, namely an open subset of $X\times X$. Our
first task is to see that this set is non-empty. This is the statement of
the following (easy and well-known) lemma. 

\proclaim Lemma 2.1. Two general lines of $X$ are skew.

{\it Proof:} If any two lines of $X$ meet, then take two of them, say
$L_1,L_2$. They meet in a point $P\in\P^{2n+1}$ and span a plane $\Pi$.
Then, from the irreducibility of $X$, either all lines of $X$ are contained
in $\Pi$ or pass through $P$. The first possibility is impossible, either
by dimensional reasons (if $n\ge3$) or by the hypothesis that $X$ is
nondegenerate. In the second possibility, $X$ will consist of the
generators of a cone with vertex $P$ over a projective
$n$-dimensional subvariety
$Y$. The nondegeneracy hypothesis implies that $Y$ spans a hyperplane in
$\P^{2n+1}$, and the fact that $X$ can be projected implies that $Y$ can be
isomorphically projected into $\P^n$. But this implies that $Y$ is a linear
space, contradicting again the nondegeneracy hypothesis.
\qed

\bigskip

\noindent {\bf Definition:} We will call {\it secant variety} to $X$ to the
variety $SX\subset G(3,2n+1)$ consisting of the closure of the set of linear
spaces spanned by pairs of skew lines of $X$. In other words, $SX$ is the
closure of the image of the rational map $p:X\times X\rationalmap G(3,2n+1)$
that associates to each pair of skew lines its linear span. We will call
the {\it secant defect} of $X$ to the dimension $\delta$ of $Y_{\Pi}$ for a
general $\Pi\in SX$, where $Y_{\Pi}$ is the set of lines of $X$ contained in
$\Pi$. It is clear, by looking at the map $p$, that $\dim(SX)=2n-2\delta$.

\noindent {\bf Observation:} These definitions are different from the
``natural'' generalization of the notion of projective secant variety and
defect. However their behavior will play a similar role as their corresponding
projective concepts, as we will see throughout the paper.

\bigskip

\proclaim Lemma 2.2. The variety $SX$ has dimension at most $2n-2$, or
equivalently, $X$ has positive secant defect.

{\it Proof:} We need to study the dominant rational map $p:X\times
X\rationalmap SX$ and show that its differential $dp_{(\l_1,\l_2)}$ at a
general point $(\l_1,\l_2)$ has rank at most $2n-2$. So assume for
contradiction that $dp$ is injective for a general
$(\l_1,\l_2)\in X\times X$. {}From Lemma 2.1 we know  that the corresponding
lines
$L_1,L_2$ are skew. We choose projective coordinates $x_0,\ldots,x_{2n+1}$
so that these two lines are 
$
L_1:x_2=\ldots=x_{2n+1}=0
$
and
$
L_2:x_0=x_1=x_4=\ldots=x_{2n+1}=0.
$
An affine chart for
$G(1,2n+1)$ around $\l_1$ consists of the lines in $\P^{2n+1}$ spanned by
the points whose coordinates are the rows of the matrix
$$
\pmatrix{1&0&a_{02}&\ldots&a_{0,2n+1}\cr
         0&1&a_{12}&\ldots&a_{1,2n+1}}
$$
(the coordinates of the chart are the $a_{ij}$'s). Similarly, an affine
chart for $G(1,2n+1)$ around $\l_2$ consists of the lines in $\P^{2n+1}$ 
spanned by the points whose coordinates are the rows of the matrix
$$
\pmatrix{b_{00}&b_{01}&1&0&b_{04}&\ldots&b_{0,2n+1}\cr
         b_{10}&b_{11}&0&1&b_{14}&\ldots&b_{1,2n+1}}
$$
Take also as an affine chart for $G(3,2n+1)$ around
$<L_1,L_2>$ to be the set of all three-spaces spanned by the rows of the matrix
$$\pmatrix{1&0&0&0&p_{04}&\ldots&p_{0,2n+2}\cr
           0&1&0&0&p_{14}&\ldots&p_{1,2n+2}\cr
           0&0&1&0&p_{24}&\ldots&p_{2,2n+2}\cr
           0&0&0&1&p_{34}&\ldots&p_{3,2n+2}}
$$
Then it is not difficult to check that the differential $dp$ at $(\l_1,\l_2)$
--the origin in our system of coordinates-- is given in these
coordinates by the equation
$$(a_{02},\ldots,a_{0,2n+1},a_{12},\ldots,a_{1,2n+1};
   b_{00},\ldots,b_{0,2n+1},b_{10},\ldots,b_{1,2n+1})\mapsto$$
$$(a_{04},\ldots,a_{0,2n+1},a_{14},\ldots,a_{1,2n+1};
   b_{04},\ldots,b_{0,2n+1},b_{14},\ldots,b_{1,2n+1})$$
The fact that $dp$ is injective is equivalent to the fact that the maps
$T_{\l_1}X\to\C^{4(n-1)}$ and $T_{\l_2}X\to\C^{4(n-1)}$ defined respectively 
by
$$(a_{02},\ldots,a_{0,2n+1},a_{12},\ldots,a_{1,2n+1})\mapsto
(a_{04},\ldots,a_{0,2n+1},a_{14},\ldots,a_{1,2n+1})$$
$$(b_{00},\ldots,b_{0,2n+1},b_{10},\ldots,b_{1,2n+1})\mapsto
(b_{04},\ldots,b_{0,2n+1},b_{14},\ldots,b_{1,2n+1})$$
are injective.

On the other hand, for general points of $X$, it is very easy to see that we 
can assume that the maps $T_{\l_1}X\to\C^{2n}$ and $T_{\l_2}X\to\C^{2n}$ 
defined respectively by
$$(a_{02},\ldots,a_{0,2n+1},a_{12},\ldots,a_{1,2n+1})\mapsto
(a_{02},\ldots,a_{0,2n+1})$$
$$(b_{00},\ldots,b_{0,2n+1},b_{10},\ldots,b_{1,2n+1})\mapsto
(b_{00},b_{01},b_{04},\ldots,b_{0,2n+1})$$
are injective. Indeed, since $X$ is uncompressed, from Remark 1.2 we know 
that the rational map $p_H:X\rationalmap H$ is generically finite over its
image for a general hyperplane $H\subset\P^{2n+1}$. If we take $\l_1$ not to
be a ramification point of $p_H$ and choose coordinates so that $H$ has
equation
$x_1=0$ then we get the wanted hypotheses for $\l_1$. We proceed similarly for
$\l_2$. 

Also (changing coordinates if necessary)
we can assume that the composition of both maps with the same linear
projection $\C^{2n}\to\C^n$ is an isomorphism. Summing up, we can assume 
from our hypothesis that the two linear maps $d_1:T_{\l_1}X\to\C^n$ and
$d_2:T_{\l_2}X\to\C^n$ defined by
$$d_1(a_{02},\ldots,a_{0,2n+1},a_{12},\ldots,a_{1,2n+1})=
(a_{0,n+2},\ldots,a_{0,2n+1})$$
$$d_2(b_{00},\ldots,b_{0,2n+1},b_{10},\ldots,b_{1,2n+1})=
(b_{0,n+2},\ldots,b_{0,2n+1})$$
are bijective.

Now we will relate
these maps with the differential of another map. More precisely, let
us consider $IX$ to be the closure of the set
$$
\{(\l_1,\l_2,\L)\in X\times X\times G(n-1,2n+1)\ | \ 
\dim<L_1,L_2>=3,\ \dim(\L\cap<L_1,L_2>)\ge1\}
$$
and let $q:IX\to G(n-1,2n+1)$ be the natural projection. Both varieties have
the same dimension $n^2+2n$ (for the dimension of $IX$ consider the projection
onto $X\times X$, whose general fibers are Schubert varieties of dimension
$n^2$). The hypothesis that $X$ can be projected to $G(1,n+1)$ means, by
condition (*), that $q$ is not surjective, equivalently that its differential
map at any point $(\l_1,\l_2,\L)\in IX$ is not injective. Let us study this
differential map at the point $(\l_1,\l_2,\L)$, where $\L$ is the linear
subspace of equations 
$
\L:x_0=x_3,\ x_1=x_{n+2}=\ldots=x_{2n+1}=0.
$
We clearly have that $(l_1,l_2,\L)$ belongs to $IX$. We take the affine chart 
for $G(n-1,2n+1)$ around $\L$ consisting of the
$(n-1)$-linear subspaces of $\P^{2n+1}$ spaned by the points whose
coordinates are the rows of the $n\times(2n+2)$-matrix
$$
\pmatrix{x_{00}&x_{01}&1&0&\ldots&0&x_{0,n+2}&\ldots&x_{0,2n+1}\cr
         1+x_{10}&x_{11}&0&1&\ldots&0&x_{1,n+2}&\ldots&x_{1,2n+1}\cr
               &\vdots& & &\ddots& &         &\vdots&          \cr
         x_{n-1,0}&x_{n-1,1}&0&0&\ldots&1&x_{n-1,n+2}&\ldots&x_{n-1,2n+1}
}
$$

Locally at the point $(\l_1,\l_2,\L)\in X\times X\times G(n-1,2n+1)$ --which
is the origin in our coordinates-- the equations for $IX$ are given by the
maximal minors of the $(n+3)\times(2n+2)$-matrix
$$
\pmatrix
{1&0&a_{02}&a_{03}&a_{04}&\ldots&a_{0,n+1}&a_{0,n+2}&\ldots&a_{0,2n+1}\cr
 0&1&a_{12}&a_{13}&a_{14}&\ldots&a_{1,n+1}&a_{1,n+2}&\ldots&a_{1,2n+1}\cr
 b_{i0}&b_{i1}&1&0&b_{i4}&\ldots&b_{i,n+1}&b_{i,n+2}&\ldots&b_{i,2n+1}\cr
 x_{00}&x_{01}&1&0&0&\ldots&0&x_{0,n+2}&\ldots&x_{0,2n+1}\cr
 1+x_{10}&x_{11}&0&1&0&\ldots&0&x_{1,n+2}&\ldots&x_{1,2n+1}\cr
         &\vdots& & & \ddots&& &   \vdots      &&    \vdots      \cr
 x_{n-1,0}&x_{n-1,1}&0&0&0&\ldots&1&x_{n-1,n+2}&\ldots&x_{n-1,2n+1}
}
$$
for $i=0,1$. In particular, the tangent space of $IX$ at the origin 
(identifying the coordinates in the tangent space with the affine coordinates)
is easily seen to be given by the equations
$$x_{0j}=b_{0j}\ \hbox{ for }\ \ j=n+2,\ldots,2n+1$$ 

\vskip -.6cm
$$x_{1j}=a_{0j}+b_{1j}\ \hbox{ for }\ \ j=n+2,\ldots,2n+1$$ 
\noindent (To see this, just take the initial terms of the maximal minors 
defined by the first $n+2$ columns of the above matrix and any of the other
columns). But now the injectivity of the above maps
$d_1$ and $d_2$ easily implies the injectivity of $dq$, so that we get a
contradiction.
\qed

\bigskip

We can improve the above result by proving that the inequality for the
dimension is in fact an equality. We prefered to separate this result in two
parts for a later discussion of both facts in relation with the use of the
uncompressedness hypothesis. The precise statement is the following.

\proclaim Lemma 2.3. For a general $\Pi\in SX$, $Y_{\Pi}$ is the curve in
$G(1,\Pi)$ consisting of the lines of one of the rulings of a smooth quadric
in $\Pi$.

{\it Proof:} We claim first that it cannot happen that any line of $Y_{\Pi}$
meets another line of $Y_{\Pi}$ (maybe infinitely close). If this happens,
$Y_{\Pi}$ would be the union of planes containing two (maybe infinitely
close) lines of $X$. {}From condition (**)  we know that the union of such
planes in $\P^{2n+1}$ has dimension at most $n+1$. As a consequence, the union
$Z$ of all the $Y_{\Pi}$'s  has dimension at most $n+1$. Since $X$ is
uncompressed, this implies that $Z$ is also the union of all lines of $X$. 

Let us prove now by induction on $k$ that the (closure of the) union of 
the spans of $k+1$ lines of $X$ is again $Z$ for any $k$ (this would give a
contradiction since there exists a value of $k$ for which the span of $k$
general lines of $X$ must be $\P^{2n+1}$). We just proved this for $k=1$. So
take now $p$ to be a point of $\P^{2n+1}$ that is in the span of $k+1$
general lines $L_1,\ldots,L_{k+1}$ of $X$ and assume $k>1$. Then $p$ is in
the span of $L_{k+1}$ and a point $p'\in <L_1,\ldots,L_k>$. By induction
hypothesis, then $p'$ is in $Z$, i.e. there exists a line $\l$ of $X$ passing
through it. But this means that $p$ is in the span of $L$ and $L_{k+1}$,
which means in  turn that $p$ is in $Z$, as wanted. This proves the claim.

We know from Lemma 2.2 that $Y_{\Pi}$ has positive dimension. It
cannot be a surface, since this would imply that any line of $Y_{\Pi}$ would
meet infinitely many others. So assume that $Y_{\Pi}$ is a curve of degree
$d$. The proof that its degree is $d=2$ will follow a very standard argument
(see for example the proof of Lemma 5.3 in \Arrondo-Sols). Take a general
line $\ell$ of $Y_{\Pi}$ and consider the Schubert variety $Z_\ell$ in
$G(1,\Pi)$ of all lines meeting
$L$. This is a quadratic cone with vertex $\ell$ (after the Pl\"ucker
embedding) and its intersection number with $Y_\Pi$ is $d$. If $Y_\Pi$ and
$Z_\ell$ are tranversal at $\ell$, the intersection multiplicity at that
point is two, and hence there are other $d-2$ lines (counted with
multiplicity) of $X$ in $Y_\Pi$ meeting $L$. If the intersection is not
transversal, the tangent line at $\ell$ of $Y_\Pi$ is a generator of the cone.
Hence there is a plane containing $\ell$ such that the intersection of $X$
with the Schubert variety of the lines contained in that plane contains a
subscheme of length at least two. 

Therefore the only possibility (after the claim) is that $d=2$ and one
easily checks also that $Y_{\Pi}$ must consist of one the rulings of a
smooth quadric.
\qed

\bigskip

This result proves that $X$ contains too many conics. For dimension $n=2$
this is the way of showing that $X$ is the Veronese surface. For general
dimension, we need to find a lot of ``special'' divisors in $X$. Fo this we
will need to generalize the above results to the span of more than two
lines of $X$. We need first to generalize a few definitions. 

\noindent {\bf Definition:} For any $k=1,\ldots,n$ let
$r_k$ be the dimension of the span of $k+1$ general lines of $X$. We define 
the {\it $k$-secant variety} to $X$ to be the subvariety $S^kX\subset
G(r_k,2n+1)$ defined as the closure of the $r_k$-linear subspaces in
$\P^{2n+1}$ spanned by $k+1$ general lines of $X$. If $\Pi$ is a general
$r_k$-space in $S^kX$ we define the set $Y_{\Pi}$ of $X$ as the set of lines
contained in $\Pi$. Of course for $k=1$ we have $S^1X=SX$. 

The wanted divisors will be the subsets $Y_{\Pi}$ for $\Pi\in S^{n-1}X$. For
this, we will need first to show that these are indeed divisors. This is the
purpose of the next result, which more generally gives the dimension of any
secant variety.

\proclaim Proposition 2.4. For any $k=1,\ldots,n$, the span $\Pi$ of $k+1$
general lines of $X$ has dimension $2k+1$, $\dim(Y_{\Pi})=k$ and
$\dim(S^kX)=(k+1)(n-k)$. 

{\it Proof:} First we observe that it is enough to show that, if $k<n$ and
$\l_1,\ldots,\l_{k+1}$ are
$k+1$ general lines of $X$ and $\Pi'=<L_1,\ldots,L_k>$, 
$\Pi=<L_1,\ldots,L_{k+1}>$, then:
$$
\dim(Y_{\Pi})\ge\dim(Y_{\Pi'})+1\eqno{(2.1)}
$$
Indeed this inequality implies that $\dim(Y_{\Pi})\ge k$ for a general $\Pi\in
S^kX$; in particular, $\dim(Y_{\Pi})\ge n$ for $\Pi\in S^nX$. Since
$X$ has dimension $n$ this shows that $Y_{\Pi}=X$ for $\Pi\in S^nX$ and we
have that all inequalities in (2.1) are equalities. This proves
$\dim(Y_{\Pi})=k$. Also, since $Y_{\Pi}=X$ for $\Pi\in S^nX$ and $X$ is
nondegenerate, it must be $\dim(\Pi)=2n+1$, from which we conclude that for
$\Pi\in S^kX$ it is $\dim(Y_{\Pi})=2k+1$. Finally, we immediately see that
the dimension of the fiber of $X\times\smash{
   \mathop{\dots}\limits^{^{k+1)}}}\times
X\rationalmap G(2k+1,2n+1)$ (the rational map assigning to $k+1$
general lines of $X$ its linear span) has dimension $k(k+1)$. This proves 
that $\dim(S^kX)=(k+1)(n-k)$. 

So let us prove inequality (2.1). Consider $\l_1,\ldots,\l_k$ to be general 
lines of $X$. Since $X$ is nondegenerate and the corresponding lines
$L_1,\ldots,L_k$ span at most a linear space $\Pi'$ of dimension
$2k-1<2n+1$ it holds that $Y_{\Pi'}$ is not the whole $X$. So we can take
a general $\l_{k+1}$ in $X$ that is not in $Y_{\Pi'}$. Consider the set $J$
to be the closure in $Y_{\Pi}\times Y_{\Pi'}$ of
$$
\{(\l,\l')\in Y_{\Pi}\times Y_{\Pi'}\ | \ L\subset<L',L_{k+1}>, \ \dim
<L',L_{k+1}>=3\}
$$
The fiber of the natural projection $J\to Y_{\Pi'}$ over a general $\l'\in
Y_{\Pi'}$ is the set $Y_{<L',L_{k+1}>}$, which has dimension one after
Lemma 2.3. Hence $\dim(J)=\dim(Y_{\Pi'})+1$. So it is enough to show that the
projection map from $J$ to
$Y_{\Pi}$ is generically finite over its image. This is so because, given a
general
$\l\in Y_{\Pi}$, a general element $\l'\in Y_{\Pi'}$
verifies that $L\subset<L',L_{k+1}>$ if $L'\subset <L,L_{k+1}>$. If the
intersection of $<L,L_{k+1}>$ with $\Pi'$ is just one line, this is precisely
$L'$. If  the intersection is a plane, $L'$ is a line in this plane that is
also in
$X$. But from Lemma 2.3, since $L$ and $L'_{k+1}$ are general, the set of
lines of
$X$ that are in $<L,L_{k+1}>$ is one of the rulings of a quadric. So there is
only one of these lines contained in a plane.
\qed

The following easy lemma can be considered as
a (partial) generalization to our context of Terracini's lemma.

\bigskip

\proclaim Lemma 2.5. Let $\Pi$ be a general element of $S^{n-1}X$. Consider
the divisor $H_{\Pi}$ of $G(1,2n+1)$ given by the Schubert cycle consisting
of all lines meeting $\Pi$. Then the intersection cycle of $X$ with $H_{\Pi}$
contains
$Y_{\Pi}$ with multiplicity at least two.

{\it Proof:} This follows immediately from the observation that the
singular locus of $H_{\Pi}$ is the Schubert cycle of all lines contained in
$\Pi$, so that $Y_{\Pi}$ is contained in that singular locus.
\qed
\bigskip
\bigskip
{\bf 3. The main theorem.}
\bigskip
We can now state and prove the main theorem of this paper.

\proclaim Theorem 3.1. Let $X$ be a smooth irreducible 
subvariety of $G(1,2n+1)$ ($n\ge 1$). Assume that $X$ is uncompressed and
nondegenerate and that it can be isomorphically projected to $G(1,n+1)$. Then
$X$  is a Veronese variety as in Example 1.1.

{\it Proof:} We keep the same notation as in the previous
sections. We will assume $n\ge 2$ since the proof for $n=1$ goes differently
and it is much easier (see \Yo). The idea is to prove that the linear
system corresponding to the hyperplanes of $\P^{2n+1}$ is the set of
divisors $Y_{\Pi}$ for $\Pi\in S^{n-1}X$ (it is a nice exercise to contrast
each step of the proof with the actual behavior of the Veronese
embedding).

Take a general $\Pi\in S^{n-1}X$. {}From Lemma 2.5 we have that
$$
H_{\Pi}\vert_X=rY_{\Pi}+E_{\Pi}\eqno{(3.1)}
$$
where $r\ge 2$ and the support of $E_{\Pi}$ does not contain $Y_{\Pi}$.

Take now a general $\Pi'\in SX$. We know from Lemma 2.3 that $Y_{\Pi'}$
consists of one of the rulings of a smooth quadric in $\Pi'$. Hence
intersecting with $Y_{\Pi'}$ in (3.1) we obtain that $2=rY_{\Pi}\cdot Y_{\Pi'}+
E_{\Pi}\cdot Y_{\Pi'}$. {}From this we see that $r=2$,
$Y_{\Pi}\cdot Y_{\Pi'}=1$ and
$E_{\Pi}\cdot Y_{\Pi'}=0$. This last equality easily implies that
$E_{\Pi}=0$. Indeed, if there exists $\l\in E_{\Pi}$, it is not difficult
to find a $\Pi'$ such that $Y_{\Pi'}$ is irreducible and is not contained in
$E_{\Pi}$; therefore the intersection of $Y_{\Pi'}$ and $E_{\Pi}$ would be
proper, hence empty since the intersection number is zero, which is a
contradiction.

So we have arrived to the equality $H_{\Pi}\vert_X=2Y_{\Pi}$. This
easily implies that all the divisors $Y_{\Pi}$ are linearly equivalent. Let us
study the complete linear system $|Y_{\Pi}|$. It clearly has no base points,
since for any point $\l\in X$ we know from Proposition 2.4 that we can find
$\l_1,\ldots,\l_n$ such that $L,L_1,\ldots,L_n$ span $\P^{2n+1}$. This
means that $l\notin Y_{\Pi}$, where $\Pi=<L_1,\ldots,L_n>$. Hence
$|Y_{\Pi}|$ defines a regular map 
$$
\varphi:X\to\P^N
$$
where $N=\dim|Y_{\Pi}|$. Let us denote by $X'$ the image of $\varphi$.
{}From what we have just seen, $|2Y_{\Pi}|$ is the hyperplane section of $X$
(after the Pl\"ucker embedding), so that $|Y_{\Pi}|$ is ample, and therefore
the map $\varphi$ is finite over $X'$. Recall that we have also got from (3.1)
the equality $Y_{\Pi}\cdot Y_{\Pi'}=1$ for general $\Pi\in S^{n-1}X$, $\Pi'\in
SX$. Hence the image of a general $Y_{\Pi'}$ is a line in $\P^N$. This proves
that two general points of $X'$ can be joined by a line, and hence
$X'=\P^n$ and $N=n$. 

On the other hand, let $E$ be the rank-two vector bundle on $X$ giving the
embedding of $X$ in $G(1,2n+1)$. Since $X$ is nondegenerate we have that
$$
m+1:=h^0(X,E)\ge 2n+2
$$ 
(Also we could conclude a priori that
equality holds, since the proof of Proposition 2.4 works in fact for $X$ in
any $G(1,m)$ with $m\ge n+3$, and we could show then that
all lines are contained in the linear span of $n+1$ general lines, hence
$m=2n+1$, from the nondegeneracy hypothesis). We also have that
$\bigwedge^2E=\O_X(2Y_{\Pi})$. Now let us show that the equality
$H_{\Pi}\vert_X=2Y_{\Pi}$ implies the splitting 
$$
E\cong\O_X(Y_{\Pi})\oplus\O_X(Y_{\Pi})\eqno{(3.2)}
$$
Indeed take a general $\Pi\in S^{n-1}X$. It is a linear space of codimension
two in $\P^{2n+1}$, hence there are two independent sections of $E$
vanishing on $Y_{\Pi}$ (therefore we can consider them as sections of
$E(-Y_{\Pi})$ ). But the equality
$H_{\Pi}\vert_X=2Y_{\Pi}$ means that the dependency locus of these two
sections is precisely $Y_{\Pi}$. In other words, there are two independet
sections of $E(-Y_{\Pi})$ whose dependency locus is empty. This proves 
$E(-Y_{\Pi})\cong\O_X\oplus\O_X$, which is (3.2).

Now (3.2) implies that there is a commutative diagram
$$
\matrix{X&\mapright{\varphi}&\P^n\cr
        \downarrow{}&\searrow&\downarrow{}\cr
        G(1,m)&\rationalmap&G(1,n+1)}
$$
Here the vertical maps are the respective embeddings of $X$ an
$\P^n$ given by the bundles $E\cong\O_X(Y_{\Pi})\oplus\O_X(Y_{\Pi})$ and
$\O_{\P^n}(1)\oplus\O_{\P^n}(1)$; the horizontal dashed arrow is the
projection induced by a linear projection $\P^m\rationalmap\P^{2n+1}$; and 
the composed diagonal morphism is the given inclusion of $X$ in $G(1,2n+1)$. 
Therefore $\varphi$ is also an embedding, hence an isomorphism, so that $X$ 
is the Veronese variety in $G(1,2n+1)$.
\qed

\bigskip
\bigskip
{\bf 4. Remarks and questions on projections of Grassmannians.}
\bigskip
One of my main goals when I started to think of this problem was to prove
Conjecture 0.1 for $n=3$. In this dimension, it follows easily from Theorem 
3.1, as we show next.

\proclaim Corollary 4.1. Let $\bar X$ be a smooth irreducible threefold of
$G(1,4)$ that is a projection of a nondegenerate threefold $X$ in $G(1,7)$. 
Then $\bar X$ is the Veronese threefold.

{\it Proof:} If $\bar X$ is not the Veronese threefold, we
know from Theorem 3.1.that $X$ must be compressed.  In other words, through a
general point of $\P^4$ there passes no line of
$\bar X$. But smooth threefolds in $G(1,4)$ verifying this property are
classified in \ABT\  and it is easy to check that none of them comes from
$G(1,7)$. 
\qed
\bigskip

In general, the natural way of approaching compressed subvarieties of
Grassmannians is the philosophy that a projective variety containing
too many lines either has a bounded degree or contains many linear varieties
of bigger dimension. This is in fact the method used in \ABT\ in dimension
three. However, for higher dimension, I do not know of any sufficiently strong
result for $n$-dimensional varieties containing an $n$-dimensional family of
lines.

There are only two places in which we used the uncompressedness hypothesis: 
in Lemmas 2.2 and 2.3. It does not seem so important its use in Lemma 2.3.
More precisely, it seems possible to prove first Proposition 2.4 without using
Lemma 2.3 (we used this lemma just for a small detail at the end of the
proof of the proposition); then one could deduce the lemma from the general
position statement in the proposition (in fact I had originally an incorrect
proof in this way, and I decided to fix the gap in the way shown in the paper
as soon as I realized that I would need uncompressedness anyway). The crucial
point where uncompressedness is used seems to be Lemma 2.2, at least the proof
strongly needs this condition. I do not know of any example of a compressed
$n$-variety projectable from $G(1,m)$ to $G(1,n+1)$ with $m\ge n+3$. If such
an example exist, it would be interesting to see whether the secant defect
is positive or not. The only example I know is the following for $m=n+2$.

\noindent {\bf Example 4.2:} Consider the embedding of $\P^1$ in $G(r,2r+2)$
given by $\O_{\P^1}(1)^{\oplus r}\oplus\O_{\P^1}(2)$, or equivalently the
smooth rational normal scroll of dimension $r+1$ in $\P^{2r+2}$. This gives a
one-dimensional family of pairwise disjoint $r$-spaces in $\P^{2r+2}$. 
Dually, we find a one-dimensional family of $(r+1)$-spaces in $\P^{2r+2}$ 
such that any two of them meet only at one point. Hence the set of lines
contained in these $(r+1)$-spaces forms a smooth $(2r+1)$-subvariety $X$ in
$G(1,2r+2)$, which is compressed if $r>0$. This variety is in fact projected
from $G(1,2r+3)$. Indeed the map $\P^1\to G(r+1,2r+2)$ correponding to the
family of $(r+1)$-spaces is defined by the epimorphism appearing in the exact
sequence
$$0\to \O_{\P^1}(-1)^{\oplus r}\oplus\O_{\P^1}(-2)\mapright{\psi}
\O_{\P^1}^{\oplus 2r+3}\to \O_{\P^1}(1)^{\oplus r+2}\to 0$$
where $\psi$ is the dual of the map defining the given embedding of $\P^1$ in
$G(r,2r+2)$. Since $h^0(\P^1,\O_{\P^1}(1)^{\oplus r+2})=2r+4$, it follows that
this family of $(r+1)$-spaces, and hence $X$, comes from $\P^{2r+3}$.

\bigskip

It could be very risky to conjecture that any compressed $n$-variety of
$G(1,n+1)$ that comes projected from a bigger $(1,m)$ is one of those in the
above example 4.2. However, having a classification of such varieties would 
certainly prove (or maybe disprove?) Conjecture 0.1. 

\bigskip
\bigskip

I would like to discuss now a little bit what should be the main items in a
general theory of secant varieties in Grassmannians of lines. Of course in 
order for the theory to work properly we need to make some general assumptions.
This includes not only the uncompressedness hypothesis (which we were unable
to avoid in our main theorem) but also it is sometimes useful to assume some
general position hypothesis. In the particular case we studied in this paper,
because of the hypothesis that our variety was nondegenerate in
$G(1,2n+1)$ we were able to prove this general position statement 
(Proposition 2.4). Let us give the precise definition of this hypothesis.

\noindent {\bf Definition:} We will say that a subvariety $X$ of $G(1,N)$ is 
in {\it general position} if for $k=1,\ldots,[{N-1\over2}]$ the span of $k+1$
general lines of $X$ is a linear space of dimension $2k+1$. Hence the $k$-th
secant variety $S^kX$ is a subvariety of $G(2k+1,N)$.

It is immediate to check that the same proof of Lemma 2.2 works to prove that 
if $N\ge n+3$, $X$ is not compressed, is not a cone and can be isomorphically
projected to 
$G(1,n+1)$ then the secant variety $SX=S^1X$ has dimension at most $2n-2$. For
general secant varieties, one can make the following definition:

\noindent {\bf Definition:} Let $X$ be a subvariety of $G(1,N)$ in general
position. We will call the {\it $k$-th secant defect} of $X$ to be the
dimension $\delta_k$ of a general $Y_{<L_1,\ldots,L_{k+1}>}$, where $Y_{\Pi}$
denotes the set of lines of $X$ contained in a given linear space $\Pi$.

Looking at the image and fibers of the rational map
$X\times\smash{
   \mathop{\dots}\limits^{^{k+1)}}}\times
X\rationalmap G(2k+1,N)$ that 
associates to each $k+1$ general lines its linear span we easily obtain that 
$$\dim(S^kX)=(k+1)(n-\delta_k)$$
Of course, if $X$ is uncompressed, then $X$ is projectable if and only if
$\delta_1>0$.

There is another relation among the defects, which is the translation to
Grassmannians of the so-called Zak's superadditivity theorem
(see \Zakpaper, or \HR, or \Fantechi). It is just the generalization of
inequality (2.1) in the proof of Proposition 2.4 and its proof is
surprisingly easy, contrary to what happens in the projective case. It is
the following.

\proclaim Proposition 4.3. If $X$ is in general position, then, for
$i+j\le[{N-1\over2}]$, one has
$$\delta_{i+j}\ge\delta_i+\delta_j$$

{\it Proof:} It is exactly the same as for (2.1), but easier since we are
already assuming general position (statement that we proved simultaneously in
Proposition 2.4). We take $\l_1,\ldots,\l_{i+j}$ to be
general lines of $X$ and consider the linear spaces
$\Pi':=<L_0,\ldots,L_i>$ and $\Pi:=<L_0,\ldots,L_{i+j}>$. Define $J$ to be the
closure in $Y_{\Pi}\times Y_{\Pi'}$ of the set
$$\{(\l,\l')\in Y_{\Pi}\times Y_{\Pi'}\ | \ L\subset <L',L_{i+1},\ldots,
L_{i+j}>, \ \ \dim<L',L_{i+1},\ldots,L_{i+j}>=2j+1\}$$
For a general $\l'\in Y_{\Pi'}$, the fiber of the natural projection
$J\to Y_{\Pi'}$ is $Y_{<L',L_{i+1},\ldots,L_{i+j}>}$. Hence $J$ has
dimension $\delta_i+\delta_j$. Finally the other projection $J\to Y_{\Pi}$
is generically finite over its image. Indeed for a general $\l$ in this
image, from the general position hypothesis we have that the spaces
$Pi'=<L_0,\ldots,L_i>$ and $<L_{i+1},\ldots,L_{i+j},L>$ meet only along a
line $L'$. This gives a unique point $(\l,\l')$ in $J$ mapping to $\l$.
\qed

\bigskip
\bigskip

{\bf References.}
\bigskip

\item{\Ad} {\AA}dslandvik, B., {\it Varieties with an extremal number of
degenerate higher secant varieties}, Journal reine angew. Math., {\bf 392}
(1987), 213-222.

\item{\Yo} Arrondo, E., {\it Subvarieties of Grassmannians}, Lecture Notes
Series  Dipartimento di Matematica Univ. Trento, {\bf 10}, 1996.

\item{\ABT} Arrondo, E. -- Bertolini, M. -- Turrini, C., {\it Congruences of
small degree in $G(1,4)$}, Preprint 1996.

\item{\Arrondo-Sols} Arrondo, E. -- Sols, I., {\it On congruences of lines in
the projective space}, M\'em. Soc. Math. France, {\bf 50}, 1992.

\item{\Fantechi} Fantechi, B., {\it On the superadditivity of secant defects},
Bull. Soc. Math. France, {\bf 118} (1990), 85-100.

\item{\HR} Holme, A. -- Roberts, J., {\it Zak's theorem on superadditivity},
Ark. Mat., {\bf 32} (1994), 99-120.

\item{\Severi} Severi, F., {\it Intorno ai punti doppi impropri di una
superficie generale dello spazio a quattro dimensioni, e a suoi punti
tripli apparenti}, Rend. Circ. Mat. Palermo, II, Ser. {\bf 15} (1901),
377-401.

\item{\Zakpaper} Zak, F.L., {\it Linear systems of hyperplane sections on 
varieties of low codimension}, Functional Anal. Appl. {\bf 19} (1986),
165-173.

\item{\Zakbook} Zak, F.L., {\it Tangents and Secants of Algebraic Varieties},
Transl. Math. Monographs AMS {\bf 127}, 1993.

\end